\documentclass[11pt]{article}
\date{10 IV 2019}

\usepackage{a4wide,latexsym,graphicx,amsmath,varioref}

\title{A note on solving nonlinear optimization problems\\ in variable precision}

\author{
  S. Gratton
  \thanks{Universit\'e de Toulouse, INP, IRIT, Toulouse, France. Email: serge.gratton@enseeiht.fr },
  and Ph. L. Toint
  \thanks{NAXYS, University of Namur, Namur, Belgium. Email: philippe.toint@unamur.be.
          Partially supported by ANR-11-LABX-0040-CIMI within the program ANR-11-IDEX-0002-02.}
}
\newcommand{\beqn}[1]{\begin{equation}\label{#1}}
\newcommand{\eeqn}{\end{equation}}
\newcommand{\req}[1]{(\ref{#1})}
\newcommand{\ms}{\;\;\;\;}
\newcommand{\tim}[1]{\;\; \mbox{#1} \;\;}

\newtheorem{theorem}{Theorem}[section]
\newtheorem{lemma}[theorem]{Lemma}
\newcommand{\numsection}[1]{\section{#1}\setcounter{equation}{0}}
\newcommand{\appnumsection}[1]{\section*{#1}\setcounter{equation}{0}
  \renewcommand{\theequation}{A.\arabic{equation}}
  \renewcommand{\thetheorem}{A.\arabic{theorem}}
  \renewcommand{\thetable}{A.\arabic{table}}
  \renewcommand{\thefigure}{A.\arabic{figure}}
  \renewcommand{\thesection}{A} }

\renewcommand{\theequation}{\arabic{section}.\arabic{equation}}
\renewcommand{\thetable}{\arabic{section}.\arabic{table}}
\renewcommand{\thefigure}{\arabic{section}.\arabic{figure}}

\newcounter{algo}[section]
\renewcommand{\thealgo}{\thesection.\arabic{algo}}
\newcommand{\algo}[3]{\refstepcounter{algo}
\begin{center}\begin{figure}[htbp]
\framebox[\textwidth]{
\parbox{0.95\textwidth} {\vspace{\topsep}
{\bf Algorithm \thealgo : #2}\label{#1}\\
\vspace*{-\topsep} \mbox{ }\\
{#3} \vspace{\topsep} }}
\end{figure}\end{center}}
\newcommand{\llem}[2]{\vspace{\baselineskip} 
\noindent\framebox[\textwidth]{\parbox{0.95\textwidth}{
\begin{lemma} \label{#1} \rm #2 \end{lemma} } } \vspace{\baselineskip} }
\newcommand{\lthm}[2]{\vspace{\baselineskip} 
\noindent\framebox[\textwidth]{\parbox{0.95\textwidth}{
\begin{theorem} \label{#1} \rm #2 \end{theorem} } } \vspace{\baselineskip} }
\newcommand{\bpr}{{\bf Proof.} \hspace{1.5mm}}
\newcommand{\epr}{\hfill $\Box$ \vspace*{1em}}

\newcommand{\proof}[1]{
\begin{list}{}{
\setlength{\topsep}{0.0pt}
\setlength{\partopsep}{0.0pt}
\setlength{\leftmargin}{0.025\textwidth}
\setlength{\rightmargin}{0.5\leftmargin}
\setlength{\labelwidth}{0.5\leftmargin}
\setlength{\labelsep}{0.25\leftmargin}}
\item \bpr #1 \epr \noindent
\end{list}}

\newcommand{\bigfrac}[2]{\frac{\displaystyle #1}{\displaystyle #2}}
\newcommand{\bigsum}{\displaystyle \sum}
\renewcommand{\Re}{\hbox{I\hskip -2pt R}}
\newcommand{\smallRe}{\hbox{\footnotesize I\hskip -2pt R}}
\newcommand{\barf}{\overline{f}}
\newcommand{\barg}{\overline{g}}
\newcommand{\sfrac}[2]{{\scriptstyle \frac{#1}{#2}}}
\newcommand{\half}{\sfrac{1}{2}}
\newcommand{\eqdef}{\stackrel{\rm def}{=}}
\newcommand{\bctable}[1]{\begin{table}[htbp]
                         \begin{center}
                         \begin{tabular}{#1} }
\newcommand{\ectable}[1]{\end{tabular}
                         \caption{#1}
                         \end{center}
                         \end{table}}
\newcommand{\tenth}{\sfrac{1}{10}}
\newcommand{\calO}{{\cal O}}
\newcommand{\calS}{{\cal S}}
\newcommand{\calU}{{\cal U}}

\newcommand{\iibe}[2]{\{ #1, \ldots, #2 \}}

\begin{document}

\maketitle

\begin{abstract}
  This short note considers an efficient variant of the trust-region algorithm
  with dynamic accuracy proposed Carter (1993) and by Conn, Gould and Toint
  (2000) as a tool for very high-performance computing, an area where it is
  critical to allow multi-precision computations for keeping the energy
  dissipation under control. Numerical experiments are presented indicating
  that the use of the considered method can bring substantial savings in
  objective function's and gradient's evaluation ``energy costs'' by
  efficiently exploiting multi-precision computations.
\end{abstract}


\section{Motivation and objectives}

Two recent evolutions in the field of scientific computing motivate the
present note.  The first is the growing importance of deep-learning methods
for artificial intelligence, and the second is the acknowledgement by computer
architects that new high-performance machines must be able to run the basic
tools of deep learning very efficiently. Because the ubiquitous mathematical
problem in deep learning is nonlinear nonconvex optimization, it is therefore
of interest to consider how to solve this problem in ways that are as
efficient as possible on new very powerful computers. As it turns out, one of
the crucial aspects in designing such machines and the algorithms that they
use is mastering energy dissipation.  Given that this dissipation is
approximately proportional to chip surface and that chip surface itself is
approximately proportional to the square of the number of binary digits
involved in the calculation \cite{GalaHoro11,PuGalaYangSchaHoro16,
High17,Mats18}, being able to solve nonlinear optimization problems with as
few digits as possible (while not loosing on final accuracy) is clearly of
interest. 

This short note's sole purpose is to show that this is possible and that
algorithms exist which achieve this goal and whose robustness significantly
exceed simple minded approaches. The focus is on unconstrained nonconvex
optimization, the most frequent case in deep learning applications.  Since the
cost of solving such problems is typically dominated by that of evaluating the
objective function (and derivatives if possible), our aim is therefore to
propose optimization methods which are able to exploit/specify varying levels
of preexisting arithmetic precision for these evaluations. Because of this
feature, optimization in this context differs from other better-studied
frameworks where the degree of accuracy may be chosen in a more continuous
way, such as adaptive techniques for optimization with noisy functions (see
\cite{ElstNeum97,Cart91b,Cart93}) or with functions the values and derivatives
of which are estimated by some (possibly dynamic) sampling process (see
\cite{XuRoosMaho17,CartSche17,BellGuriMori18,ChenJianLinZhan18,
  BlanCartMeniSche16,BergDiouKungRoye18}, for instance). We propose here a
suitable adaptation of the dynamic-accuracy trust-region framework proposed by
Carter in \cite{Cart93} and by Conn, Gould and Toint in Section~10.6 of
\cite{ConnGoulToin00} to the context of multi-precision computations. Our
proposal complements that of \cite{Leyfetal16}, where inexactness is also used
for energy saving purposes, and where its exploitation is restricted to the
inner linear algebra work of the solution algorithm, while still assuming
exact values of the nonlinear function involved\footnote{The solution of
  nonlinear systems of equations is considered rather than unconstrained
  optimization.}. Note that the framework of inexact computation has already
been discussed in other contexts \cite{Baboetal09,Pale14,Kugl15,Leyfetal16}.

The paper is organized as follows.  Section~\ref{algo-s}
presents the algorithmic framework using variable accuracy.
Section~\ref{experiments-s} reports encouraging numerical results suggesting
the potential of the approach, while conclusions and perspectives for further
research are discussed in Section~\ref{conclusions-s}.

\numsection{Nonconvex Optimization with Dynamic Accuracy}
\label{algo-s}

We start by briefly recalling the context of the dynamic-accuracy trust-region
technique of \cite{ConnGoulToin00}. Consider the unconstrained minimization
problem
\beqn{problem}
\min_{x \in \smallRe^n} f(x)
\eeqn
where $f$ is a sufficiently smooth function from $\Re^n$ into $\Re$, and where
the value of the objective can be approximated with a
prespecified level of accuracy. This is to say, given $x \in \Re^n$
and a an accuracy levels $\omega_f > 0$, the function $\barf(x,\omega_f)$ such that
\beqn{barf-def} 
|\barf(x, \omega_f) - \barf(x,0)| \leq \omega_f
\tim{ and }
\barf(x,0) = f(x)
\eeqn
The crucial difference with a more standard approach for optimization with
noisy functions is that \emph{the required accuracy level $\omega_f$ may be
specified by the minimization algorithm itself within a prespecified set},
with the understanding that the more accurate the requirement specified by
$\omega_f$, the higher the ``cost'' of evaluating $\barf$.

We propose to use a trust-region method, that is an iterative algorithm where,
at iteration $k$, a first-order model $m(x_k,s)$ is approximately minimized on
$s$ in a ball (the trust region) centered at the current iterate $x_k$ and of
radius $\Delta_k$ (the trust-region radius), yielding a trial point. The value
of the reduction in the objective function achieved at the trial point is then
compared to that predicted by the model.  If the agreement is deemed
sufficient, the trial point is accepted as the next iterate $x_{k+1}$ and the
radius kept the same or increased. Otherwise the trial point is rejected and
the radius reduced.  This basic framework, whose excellent convergence
properties and outstanding numerical performance are well-known, was modified
in Section~10.6 of \cite{ConnGoulToin00} to handle the situation when only
$\barf$ is known, rather than $f$.  It has already been adapted to other
contexts (see \cite{BellGratRicc18} for instance).

However, the method has the serious drawback of requiring exact gradient
values, even if function value may be computed inexactly. We may then call on
Section~8.4 of \cite{ConnGoulToin00} which indicates that convergence of
trust-region methods to first-order critical points may be ensured with
inexact gradients.  If we now define the approximate gradient at $x$ as
the function $\barg(x,\omega_{g})$ such that
\beqn{barg-def}
\|\barg(x,\omega_g)-\barg(x,0)\| \leq \omega_g \|\barg(x,\omega_g)\|
\tim{and}
\barg(x,0) = \nabla_x^1f(x),
\eeqn
this convergence is ensured provided the relative error on the gradient
remains suitably small througout the iterations,
that is
\beqn{omega-g-cond}
0 \leq \omega_g \leq \kappa_g,
\eeqn
where $\kappa_g$ is specified below in Algorithm~~\ref{TR1DA}.  Note
that this relative error need not tend to zero, but that the absolute error
will when convergence occurs (see \cite[Theorem~8.4.1,
  p.~281]{ConnGoulToin00}). Also note that the concept of a relative gradient
error is quite natural if one assumes that $\barg(x,\omega_{g})$ is computed
using an arithmetic with a limited number of significant digits.

We now propose a variant of this scheme, which we state as
Algorithm~\ref{TR1DA} \vpageref{TR1DA}.

\algo{TR1DA}{Trust region with dynamic accuracy on $f$ and $g$ (TR{\bf 1}DA)}{
  \begin{description}
  \item[Step 0: Initialization: ] An initial point $x_0$, an initial trust-region
    radius $\Delta_0$, an initial accuracy levels $\omega_{f,0}$ and a desired
    final gradient accuracy $\epsilon \in (0,1]$ are given. The positive constants
    $\eta_0$, $\eta_1$, $\eta_2$, $\gamma_1$, $\gamma_2$, $\gamma_3$ and $\kappa_g$ are also given and 
    satisfy
    \beqn{alg-consts}
    0<\eta_1\leq\eta_2<1,
    \;\;
    0<\gamma_1\leq \gamma_2 <1 \leq \gamma_3,
    \;\;
    0< \eta_0 < \half \eta_1
    \tim{ and }
    \eta_0+\kappa_g < \half(1-\eta_2).
    \eeqn
    Compute $f_0= \barf(x_0,\omega_{f,0})$ and set $k=0$.
  \item[Step 1: Check for termination: ] If $k=0$ or $x_k\neq x_{k-1}$,
    choose $\omega_{g,k} \in (0,\kappa_g]$  and compute
    $\barg(x_k,\omega_{g,k})$ such that
    \beqn{bargk-def}
    \|\barg(x_k,\omega_{g,k})-\barg(x_k,0)\| \leq \omega_{g,k} \|\barg(x_k,\omega_{g,k})\|.
    \eeqn
    In all cases, terminate if
    \beqn{term-g}
    \|\barg(x_k,\omega_{g,k})\| \leq \frac{\epsilon}{1+\kappa_g}.
    \eeqn
  \item[Step 2: Step calculation: ] Select a symmetric Hessian approximation
    $H_k$ and compute a step $s_k$ such that $\|s_k\|\leq \Delta_k$  which
    sufficiently reduces the model 
    \beqn{model}
    m(x_k,s) = f_k + \barg(x_k,\omega_{g,k})^Ts + \half s^TH_ks
    \eeqn
    in the sense that
    \beqn{cauchy}
    m(x_k,0)-m(x_k,s_k)
    \geq \half \|\barg(x_k,\omega_{g,k})\| \min\left[\frac{\|\barg(x_k,\omega_{g,k})\|}{1+\|H_k\|},\Delta_k\right]
    \eeqn
  \item[Step 3: Evaluate the objective function:]
    Select
    \beqn{omega-cond}
    \omega_{f,k}^+ \in \Big(\,0, \eta_0[m(x_k,0)-m(x_k,s_k)]\,\Big]
    \eeqn
    and compute $f_k^+ = \barf(x_k+s_k,\omega_{f,k}^+)$.
    If $\omega_{f,k}^+ < \omega_{f,k}^{\ }$,  recompute  $f_k=\barf(x_k,\omega_{f,k}^+)$ . 
  \item[Step 4: Acceptance of the trial point: ] Define the ratio
    \beqn{rho-def}
    \rho_k = \frac{f_k-f_k^+}{m(x_k,0)-m(x_k,s_k)}.
    \eeqn
    If $\rho_k \geq \eta_1$, then define $x_{k+1} = x_k+s_k$, $f_{k+1}=f_k^+$ and set
    $\omega_{f,k+1}= \omega_{f,k}^+$.  Otherwise set $x_{k+1}=x_k$,
    $\omega_{f,k+1}= \omega_{f,k}$ and $\omega_{g,k+1}= \omega_{g,k}$.
  \item[Step 5: Radius update:] Set
    \beqn{rad-update}
    \Delta_{k+1} \in \left\{ \begin{array}{ll}
      {}[\Delta_k,\gamma_3\Delta_k) &\tim{if} \rho_k \geq \eta_2, \\
      {}[\gamma_2\Delta_k, \Delta_k) & \tim{if} \rho_k \in [\eta_1,\eta_2),\\
      {}[\gamma_1\Delta_k,\gamma_2\Delta_k] & \tim{if}  \rho_k < \eta_1.
    \end{array} \right.
    \eeqn
    Increment $k$ by 1 and go to Step 1.
  \end{description}
}

This algorithm differs from that presented in \cite[p.~402]{ConnGoulToin00} on
two accounts. First, it incorporates inexact gradients, as we discussed above.
Second, it does not require that the step $s_k$ is recomputed whenever a more
accurate objective function's value $f_k=\barf(x_k,\omega_{f,k}^+)$ is
required in Step~3. This last feature makes the algorithm more efficient.
Moreover, it does not affect the sequence of iterates since the value of the
model decrease predicted by the step is independent of the objective function
value.  As a consequence, the convergence to first-order points studied in
Section~10.6.1 of \cite{ConnGoulToin00} (under the assumption that the
approximate Hessians $H_k$ remain bounded) still applies.  In what follows, we
choose to contruct this approximation using a limited-memory symmetric
rank-one (SR1) quasi-Newton update\footnote{Numerical experiments not reported
here suggest that our default choice of remembering 15 secant pairs gives
good performance, although keeping a smaller number of pairs is still
acceptable.} based on gradient differences \cite[Section~8.2]{NoceWrig99}.
Also note that condition \req{cauchy} enforces the standard ``Cauchy
decrease'' which is easily obtained by minimizing the model \req{model} in the intersection of
the trust region $\{s \in \Re^n \mid \|s\| \leq \Delta_k\}$ and the direction of the negative gradient
(see \cite[Theorem~6.3.1]{ConnGoulToin00}).

As it turns out, this variant of \cite{ConnGoulToin00} is quite close to the
method proposed by Carter in \cite{Cart93}, the main difference being that the
latter uses fixed tolerances in slighly different ranges\footnote{
  Carter \cite{Cart93} requires $\omega_g \leq 1-\eta_2$ while we require
  $\omega_g\leq \kappa_g$ with $\kappa_g$ satisfying \req{alg-consts}.
  A fixed value is also used for $\omega_f$,
  whose upper bound depends on $\omega_g$.  The Hessian approximation is
  computed using an unsafeguarded standard BFGS update.}.

We immediately note that, at termination,
\beqn{term-grad}
\|\nabla_x^1 f(x_k)\|
\leq \|\barg(x_k,\omega_{g,k})\| + \|\barg(x_k,\omega_{g,k})-\barg(x_k,0)\|
\leq (1+\omega_{g,k}) \|\barg(x_k,\omega_{g,k})\|
\leq \epsilon.
\eeqn
where we have used the triangle inequality, \req{barg-def} and
\req{omega-g-cond}. As a consequence, the TR$1$DA algorithm terminates at a
true $\epsilon$-approximate first-order-necessary minimizer.  Moreover, the
arguments leading to \cite[Theorem~8.4.5]{ConnGoulToin00} and the development
of p.~404-405 in the same reference can be combined (see the Appendix for
details) to prove that the maximum number of iterations needed by Algorithm
TR$1$DA to find such an $\epsilon$-approximate first-order-necessary minimizer
is $O(\epsilon^{-2})$.  Moreover, this bound was proved sharp in most cases in
\cite{CartGoulToin18a}, even when exact evaluations of $f$ and $\nabla_x^1f$
are used.

\numsection{Numerical Experience}\label{experiments-s}

We now present some numerical evidence that the TR$1$DA algorithm can perform
well and provide significant savings in terms of energy costs, when these are
dominated by the function and gradient evaluations. Our experiments
are run in Matlab (64bits) and use a collection of 86 small unconstrained test
problems\footnote{The collection of \cite{Buck89} and a few other problems,
all available in Matlab.} detailed in Table~\ref{testprobs-t}.
In what follows, we declare a run successful when an iterate is found such
that \req{term-g}, and hence \req{term-grad}, hold in at most 1000 iterations.

{\small 
\bctable{|lrc|lrc|lrc|}
\hline
Name & dim. & source & Name & dim. & source & Name & dim. & source \\
\hline
{\tt argauss}   &  3 & \cite{MoreGarbHill81,Buck89} &
{\tt arglina}   & 10 & \cite{MoreGarbHill81,Buck89} &
{\tt arglinb}   & 10 & \cite{MoreGarbHill81,Buck89} \\
{\tt arglinc}   & 10 & \cite{MoreGarbHill81,Buck89} &
{\tt argtrig}   & 10 & \cite{MoreGarbHill81,Buck89} &
{\tt arwhead}   & 10 & \cite{ConnGoulLescToin94,GoulOrbaToin15b}\\
{\tt bard}      &  3 & \cite{Buck89} &
{\tt bdarwhd}   & 10 & \cite{GoulOrbaToin15b} &
{\tt beale}     &  2 & \cite{Buck89}\\
{\tt biggs6}    &  6 & \cite{MoreGarbHill81,Buck89} &
{\tt box}       &  3 & \cite{MoreGarbHill81,Buck89} &
{\tt booth}     &  2 & \cite{Buck89} \\
{\tt brkmcc}    &  2 & \cite{Buck89} &
{\tt brazthing} &  2 & - &
{\tt brownal}   & 10 & \cite{MoreGarbHill81,Buck89} \\
{\tt brownbs}   &  2 & \cite{MoreGarbHill81,Buck89} &
{\tt brownden}  &  4 & \cite{MoreGarbHill81,Buck89} &
{\tt broyden3d} & 10 & \cite{MoreGarbHill81,Buck89} \\
{\tt broydenbd} & 10 & \cite{MoreGarbHill81,Buck89} &
{\tt chebyqad}  & 10 & \cite{MoreGarbHill81,Buck89} &
{\tt cliff}     &  2 & \cite{Buck89} \\
{\tt clustr}    &  2 & \cite{Buck89} &
{\tt cosine}    &  2 & \cite{GoulOrbaToin15b} &
{\tt crglvy}    & 10 & \cite{MoreGarbHill81,Buck89} \\
{\tt cube}      &  2 & \cite{Buck89} &
{\tt dixmaana}  & 12 & \cite{DixoMaan88,Buck89} &
{\tt dixmaanj}  & 12 & \cite{DixoMaan88,Buck89} \\
{\tt dixon}     & 10 & \cite{Buck89} &
{\tt dqrtic}    & 10 & \cite{Buck89} &
{\tt edensch}   &  5 & \cite{Li88} \\
{\tt eg2}       & 10 & \cite{ConnGoulToin92,GoulOrbaToin15b} &
{\tt eg2s}      & 10 & \cite{ConnGoulToin92,GoulOrbaToin15b} &
{\tt engval1}   & 10 & \cite{MoreGarbHill81,Buck89} \\
{\tt engval2}   & 10 & \cite{MoreGarbHill81,Buck89} &
{\tt freuroth}  &  4 & \cite{MoreGarbHill81,Buck89} &
{\tt genhumps}  &  2 & \cite{GoulOrbaToin15b}\\
{\tt gottfr}    &  2 & \cite{Buck89} &
{\tt gulf}      &  4 & \cite{MoreGarbHill81,Buck89} &
{\tt hairy}     &  2 & \cite{GoulOrbaToin15b}\\
{\tt helix}     &  3 & \cite{MoreGarbHill81,Buck89} &
{\tt hilbert}   & 10 & \cite{Buck89} &
{\tt himln3}    & 10 & \cite{Buck89} \\
{\tt himm25}    & 10 & \cite{Buck89} &
{\tt himm27}    & 10 & \cite{Buck89} &
{\tt himm28}    & 10 & \cite{Buck89} \\
{\tt himm29}    & 10 & \cite{Buck89} &
{\tt himm30}    & 10 & \cite{Buck89} &
{\tt himm33}    & 10 & \cite{Buck89} \\
{\tt hypcir}    &  2 & \cite{Buck89} &
{\tt indef}     &  5 & \cite{GoulOrbaToin15b} &
{\tt integreq}  &  2 & \cite{MoreGarbHill81,Buck89} \\
{\tt jensmp}    &  2 & \cite{MoreGarbHill81,Buck89} &
{\tt kowosb}    &  4 & \cite{MoreGarbHill81,Buck89} &
{\tt lminsurf}  & 25 & \cite{GrieToin82b,Buck89}\\
{\tt mancino}   & 10 & \cite{Sped75,Buck89} &
{\tt mexhat}    &  2 & \cite{BrowBart87} &
{\tt meyer3}    &  3 & \cite{MoreGarbHill81,Buck89} \\
{\tt morebv}    & 12 & \cite{MoreGarbHill81,Buck89} &
{\tt msqrtals}  & 16 & \cite{Buck89} &
{\tt msqrtbls}  & 16 & \cite{Buck89} \\
{\tt nlminsurf} & 25 & \cite{GrieToin82b,Buck89}&
{\tt osbornea}  &  5 & \cite{MoreGarbHill81,Buck89} &
{\tt osborneb}  & 11 & \cite{MoreGarbHill81,Buck89} \\
{\tt penalty1}  & 10 & \cite{MoreGarbHill81,Buck89} &
{\tt penalty2}  & 10 & \cite{MoreGarbHill81,Buck89} &
{\tt powellbs}  &  2 & \cite{MoreGarbHill81,Buck89} \\
{\tt powellsg}  &  4 & \cite{MoreGarbHill81,Buck89} &
{\tt powellsq}  &  2 & \cite{Buck89} &
{\tt powr}      & 10 & \cite{Buck89} \\
{\tt recipe}    &  2 & \cite{Buck89} &
{\tt rosenbr}   &  2 & \cite{MoreGarbHill81,Buck89} &
{\tt schmvett}  &  3 & \cite{SchmVett70,Buck89} \\
{\tt scosine}   &  2 & \cite{GoulOrbaToin15b} &
{\tt sisser}    &  2 & \cite{Buck89} &
{\tt spmsqrt}   & 10 & \cite{Buck89} \\
{\tt tquartic}  & 10 & \cite{Buck89} &
{\tt tridia}    & 10 & \cite{Buck89} &
{\tt trigger}   &  7 & \cite{PoenSchw81,Buck89} \\
{\tt vardim}    & 10 & \cite{MoreGarbHill81,Buck89} &
{\tt watson}    & 12 & \cite{MoreGarbHill81,Buck89} &
{\tt wmsqrtals} & 16 & - \\
{\tt wmsqrtbls} & 16 & - &
{\tt woods}     & 12 & \cite{MoreGarbHill81,Buck89} &
{\tt zangwil2}  &  2 & \cite{Buck89} \\
{\tt zangwil3}  &  3 & \cite{Buck89} &
& & &
& & \\
\hline
\ectable{\label{testprobs-t}The test problems}
}

In the following set of experiments with the TR$1$DA variants,
we assume that the objective function's value $\barf(x_k,\omega_{f,k}))$ and the
gradient $\barg(x_k,\omega_{g,k})$ can be computed in double, single or half
precision (with corresponding accuracy level equal to machine precision, half
machine precision or quarter machine precision). In our experiments, single
and half precision are simulated by adding a uniformly distributed random
numerical perturbation in the ranges $[-10^{-8},10^{-8}]$ and
$[-10^{-4},10^{-4}]$, respectively. Thus, when the TR$1$DA
algorithm specifies an accuracy level $\omega_{f,k}$ or $\omega_{g,k}$, this may not be attainable as
such, but the lower of the three available levels of accuracy is then chosen to
perform the computation in (possibly moderately) higher accuracy than
requested. The \emph{equivalent double-precision costs} of the
evaluations of $f$ and $g$ in single precision are then computed by dividing
the cost of evaluation in double precision by four\footnote{Remember it is
proportional to the square of the number of significant digits.}.  Those for
half precision are correspond to double-precision costs divided by sixteen.

To set the stage, our first experiment starts by comparing three variants of the TR$1$DA
algorithm:
\begin{description}
\item[LMQN: ] a version using $\omega_{f,k} = \omega_{g,k} = 0$ for all $k$ (i.e.\ 
  using the full double precision arithmetic throughout),
\item[LMQN-s: ] a version using single precision evaluation of the
  objective function and gradient for all $k$,
\item[LMQN-h: ] a version using half precision evaluation of the
  objective function and gradient for all $k$.
\end{description}
These variants correspond to a simple minded approach where the expensive
parts of the optimization calculation are conducted in reduced precision
without any further adaptive accuracy management.
For each variant, we report, for three different values of the final gradient
accuracy $\epsilon$,
\begin{enumerate}
  \vspace*{-1mm}
\item the robustness as given by the number of successful solves for the
  relevant $\epsilon$ (nsucc),
  \vspace*{-1mm}
\item the average number of iterations (its.),
  \vspace*{-1mm}
\item the average equivalent double-precision costs of objective function's
evaluations (costf), 
\item \vspace*{-1mm}
  the average equivalent double-precision costs of gradient's
  evaluations (costg),
\item \vspace*{-1mm} the ratio of the average number of iterations used by the variant
  compared to that used by LMQN, computed on problems solved by both LMQN and
  the variant (rel. its.), 
\item   \vspace*{-1mm}
  the ratio of the average equivalent double-precision evaluation costs for $f$ used by the variant
  compared to that used by LMQN, computed on problems solved by both LMQN and
  the variant (rel. costf), 
\item   \vspace*{-1mm}
the ratio of the average equivalent double-precision evaluation costs for $g$ used by the variant
  compared to that used by LMQN, computed on problems solved by both LMQN and
  the variant (rel. costg), 
  \vspace*{-1mm}
\end{enumerate}
where all averages are computed on a sample of 20 independent runs. We are
interested in making the values in the last two indicators as small as possible
while maintaining a reasonable robustness (reported by nsucc).

\bctable{clrrrrccc}
\hline
$\epsilon$ & Variant  & nsucc & its.~~ & costf & costg & rel. its. & rel. costf & rel. costg\\
\hline
 1e-03  & LMQN   & 82 &  41.05 &   42.04 &   42.04 & & & \\ 
        & LMQN-s & 78 &  41.40 &   42.60 &   42.60 &   1.03 &    1.04 &    1.04 \\ 
        & LMQN-h & 22 &  16.95 &    1.12 &    1.12 &   0.97 &    0.06 &    0.06 \\ 
 1e-05  & LMQN   & 80 &  46.34 &   47.38 &   47.38 & & & \\ 
        & LMQN-s & 48 &  47.79 &   48.96 &   48.96 &   1.08 &    1.08 &    1.08 \\ 
        & LMQN-h & 10 &  17.80 &    1.18 &    1.18 &   1.38 &    0.08 &    0.08 \\ 
 1e-07  & LMQN   & 67 &  62.76 &   63.85 &   63.85 & & & \\ 
        & LMQN-s & 25 &  28.28 &   28.96 &   28.96 &   0.82 &    0.81 &    0.81 \\ 
        & LMQN-h &  6 &  15.83 &    1.05 &    1.05 &   0.97 &    0.06 &    0.06 \\ 
\hline
\ectable{\label{dsh-t}Results for LMQN-s and LMQN-h compared to LMQN}

Table~\ref{dsh-t} shows that the variants LMQN-s and LMQN-h compare very
poorly to LMQN for two reasons. The first and most important is the quickly
decreasing robustness when the final gradient accuracy $\epsilon$ gets
tighter. The second is that, even for the cases where the robustness is not
too bad (LMQN-s for $\epsilon= 10^{-3}$ and maybe $10^{-5}$), we observe no
improvement in costf and costg (as reported in the two last columns of the table).
However, and as expected, when LMQN-h happens to succeed, it does so at a much
lower cost, both for $f$ and $g$.

Our second experiment again compares 3 variants:
\begin{description}
\item[ LMQN:   ] as above,
\item[ iLMQN-a: ] a variant of the TR$1$DA algorithm where, for each $k$
  \beqn{omegaf-1}
  \omega_{f,k}^+ = \min[ \tenth, \sfrac{4}{100}\eta_1 \big(m(x_k,0)-m(x_k,s_k)\big) ]
  \tim{ and }
  \omega_{g,k} = \half\kappa_g.
  \eeqn
\item[ iLMQN-b: ] a variant of the TR$1$DA algorithm where, for each $k$,
  $\omega_{f,k}$ is chosen as in \req{omegaf-1} and 
  \beqn{omegag-2}
  \omega_{g,k} = \min [ \kappa_g, \omega_{f,k}].
  \eeqn
\end{description}
The updating formulae for iLMQNa are directly inspired by \req{omega-cond} and
\req{omega-g-cond} above.  The difference between the two updates of
$\omega_{g,k}$ appears to give contrasted but interesting outcomes, as we
discuss below.  The results obtained for these variants are presented in
Table~\ref{vars-t} in the same format as that used for Table~\ref{dsh-t}, the
comparison in the last three columns being again computed with respect to
LMQN.

\bctable{clrrrrccc}
\hline
$\epsilon$ & Variant  & nsucc & its.~~ & costf & costg & rel. its. & rel. costf & rel. costg\\
\hline
 1e-03  & LMQN    & 82 &  41.05 &   42.04 &   42.04 & & & \\ 
        & iLMQN-a & 80 &  50.05 &    9.88 &    6.11 &   1.23 &    0.24 &    0.15 \\ 
        & iLMQN-b & 76 &  52.67 &   13.85 &    3.34 &   1.36 &    0.35 &    0.08 \\ 
 1e-05  & LMQN    & 80 &  46.34 &   47.38 &   47.38 & & & \\ 
        & iLMQN-a & 75 &  75.92 &   36.21 &   24.77 &   1.40 &    0.63 &    0.42 \\ 
        & iLMQN-b & 63 &  72.57 &   39.85 &    4.60 &   1.78 &    0.95 &    0.11 \\ 
 1e-07  & LMQN    & 67 &  62.76 &   63.85 &   63.85 & & & \\ 
        & iLMQN-a & 47 &  65.83 &   58.97 &   37.50 &   1.18 &    1.03 &    0.65 \\ 
        & iLMQN-b & 40 &  87.35 &   95.09 &    5.52 &   1.39 &    1.45 &    0.09 \\ 
\hline
\ectable{\label{vars-t}Results for the variable-precision variants}

\begin{figure}[htbp]
\centerline{
\includegraphics[width=6.5cm]{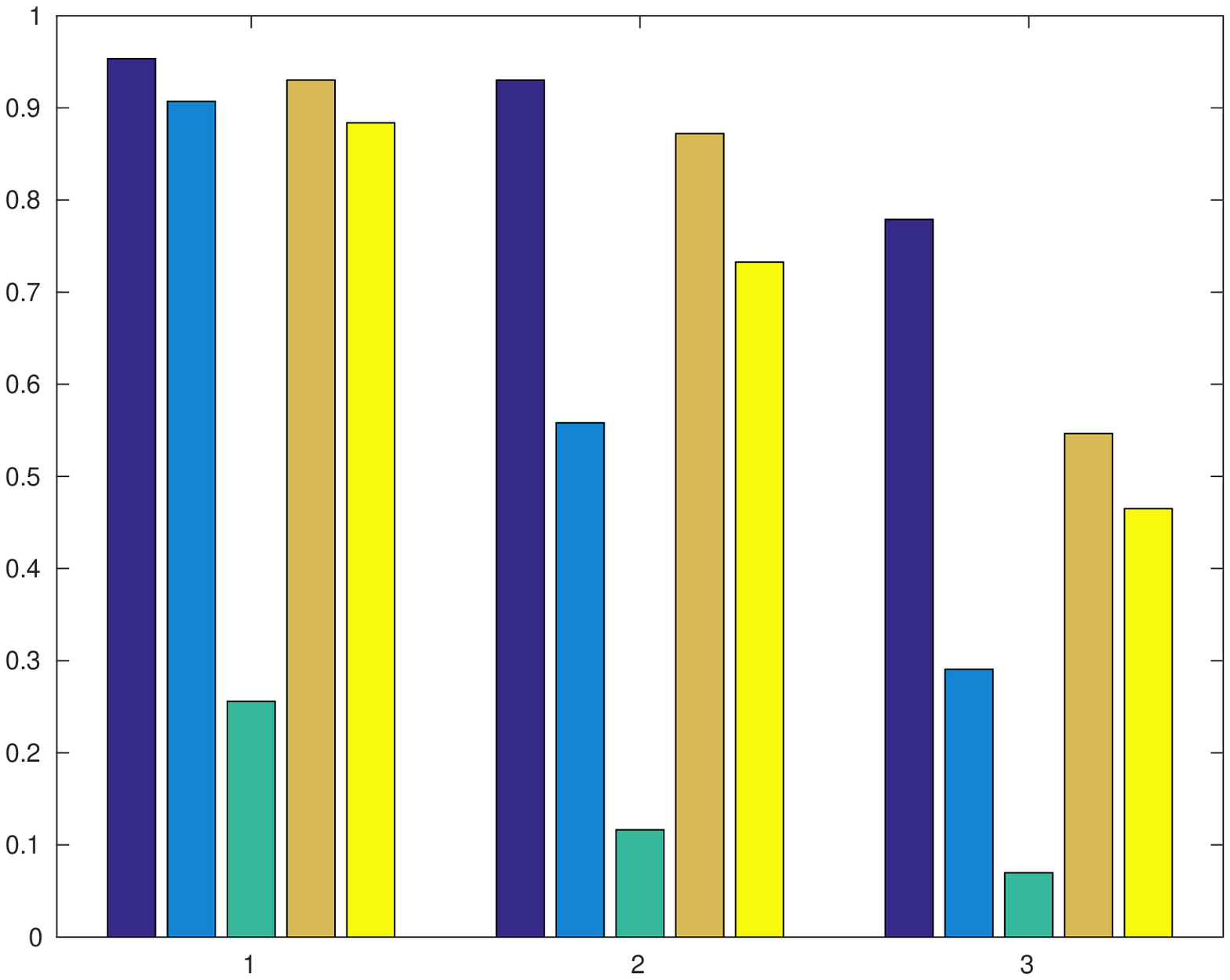}
\hspace*{1cm}
\includegraphics[width=6.5cm]{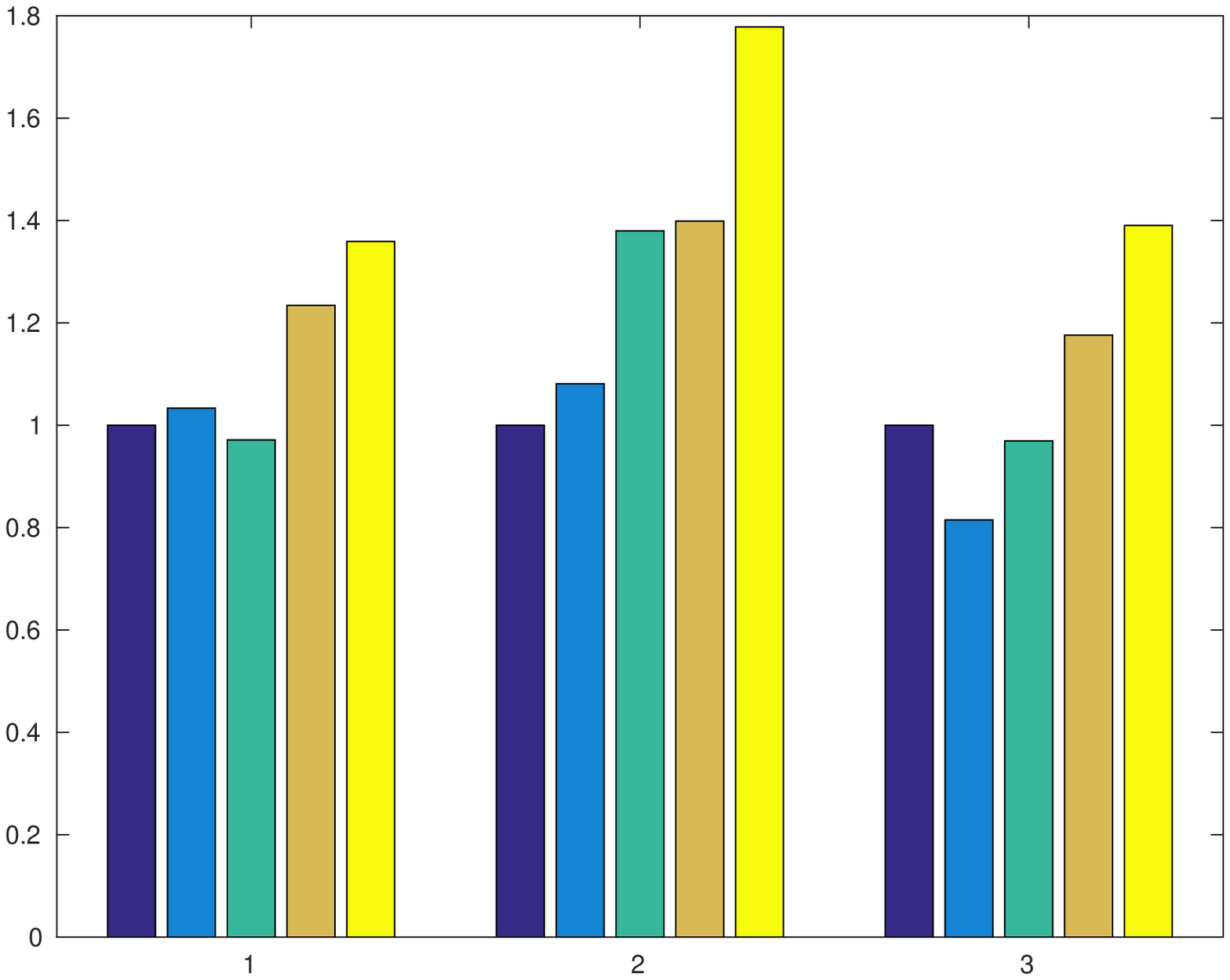}
}
\caption{\label{f_succ-its} Relative reliabilities and iteration numbers
}
\end{figure}

\begin{figure}[htbp]
\centerline{
\includegraphics[width=6.5cm]{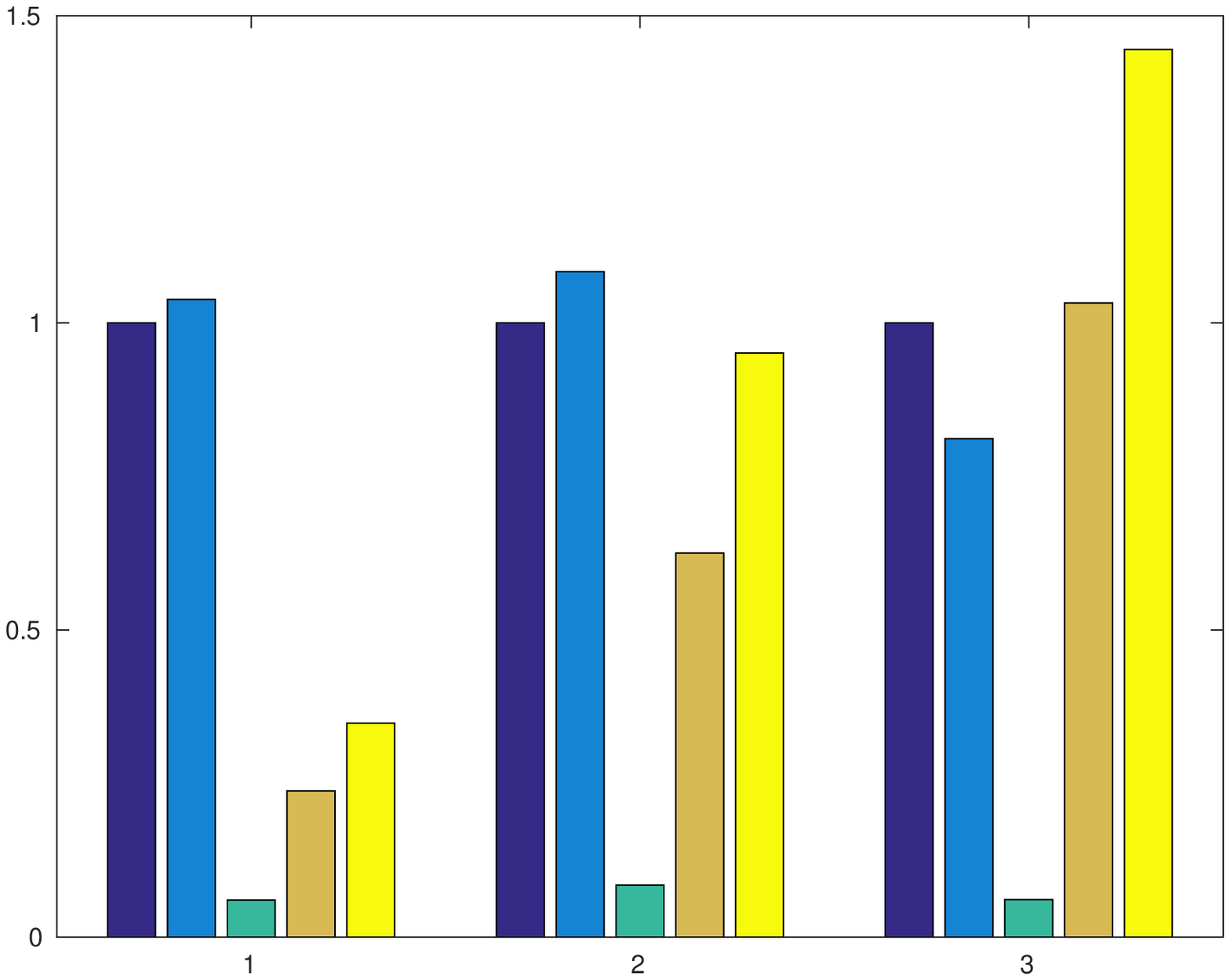}
\hspace*{1cm}
\includegraphics[width=6.5cm]{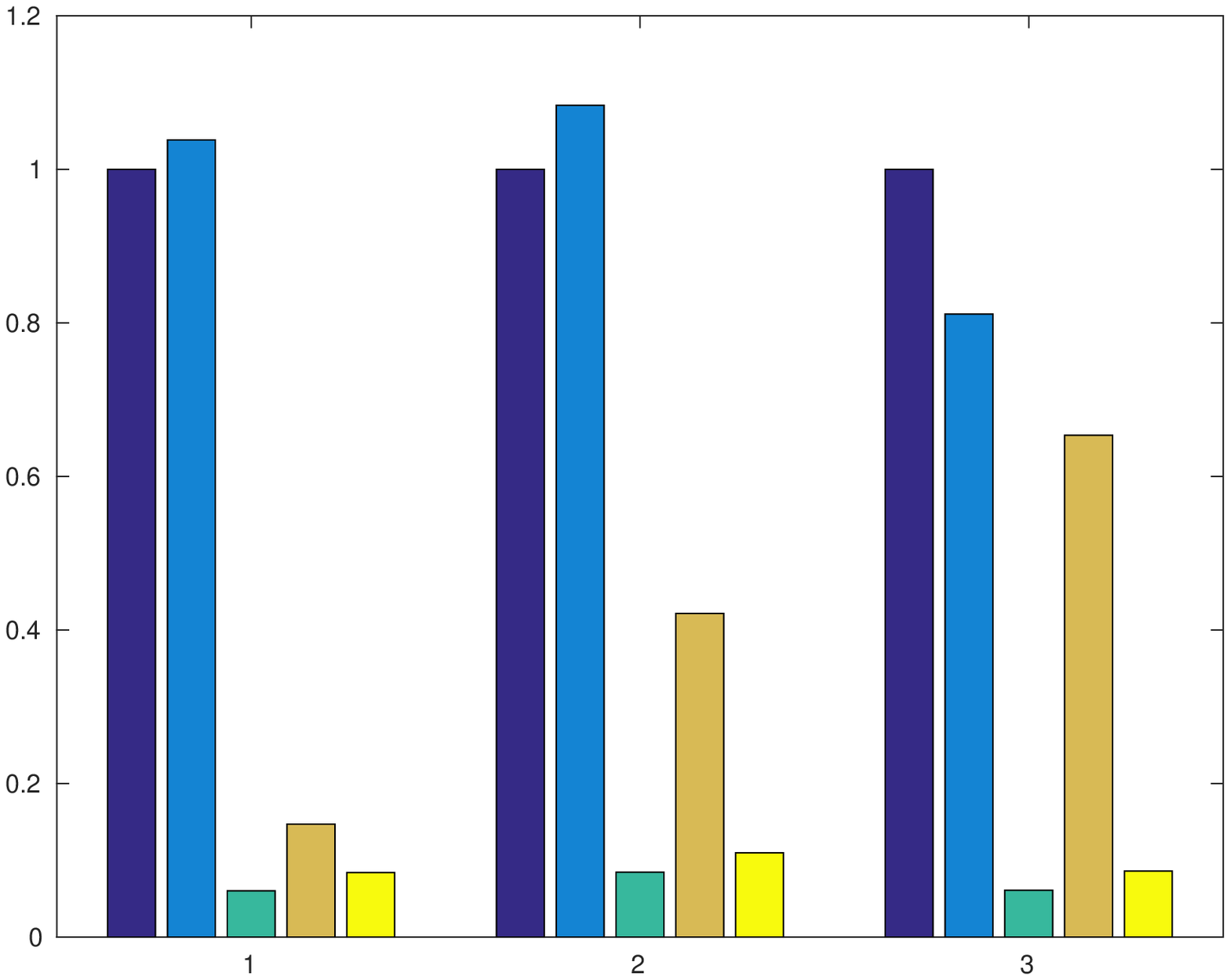}
}
\caption{\label{f_cf-cg} Relative energy savings for the evaluations of $f$
  and $g$
}
\end{figure}

These results are graphically summarized in Figures~\ref{f_succ-its} and \ref{f_cf-cg}.
In both figures, each group of vars represents the performance of the five
methods discussed above: LMQN (dark blue), LMQN-s (light blue), LMQN-h
(green), iLMQN-a (brown) and iLMQN-b (yellow).  The left part of
Figure~\ref{f_succ-its} gives the ratio of successful solves to the number of
problems, while the right part shows the relative number of iterations
compared to that used by LMQN (on problems solved by both algorithms).
Figure~\ref{f_cf-cg} gives the relative energy costs compared to LMQN, the
right part presenting the costs of evaluating the objective function and the
left part the costs of evaluating the gradients.


The following conclusions follow from Table~\ref{vars-t} and
Figures~\ref{f_succ-its}-\ref{f_cf-cg}.
\begin{enumerate}
\item For moderate final accuracy requirements ($\epsilon = 10^{-3}$ or
  $10^{-5}$), the inexact variants iLMQN-a and iLMQN-b perform well: they
  provide very acceptable robustness compared to the exact method and, most
  importantly here, yield very significant savings in costs, both for the
  gradient and the objective function, at the price of a reasonable increase in
  the number of iterations.
\item The iLMQN-a variant appears to dominate the iLMQN-b in
  robustness and savings in the evaluation of the objective function. iLMQN-b
  nevertheless shows significantly larger savings in the gradient's
  evaluation costs, but worse performance for the evaluation of the objective
  function.  
\item When the final accuracy is thigher ($\epsilon= 10^{-7}$), the inexact
  methods seem to loose their edge.  Not only they become less robust
  (especially iLMQN-b), but the gains in function evaluation costs disappear
  (while those in gradient evaluation costs remain significant for the
  problems the methods are able to solve). A closer examination of the
  detailed numerical results indicates that, unsurprisingly, inexact methods
  mostly fail on ill-conditioned problems (e.g.\ {\tt brownbs, powellbs, meyer3,
    osborneb}).
\item The comparison of iLMQN-a and even iLMQN-b with LMQN-s and LMQN-h
  clearly favours the new methods both in robustness and gains obtained,
  showing that purpose-designed algorithms outperform simple-minded approaches
  in this context.
\end{enumerate}
Summarizing, the iLMQN-a multi-precision algorithm appears, in our
experiments, to produce significant savings in function's and gradient's evaluation
costs when the final accuracy requirement and/or the problem conditioning is
moderate. Using the iLMQN-b variant may produce, on the problems where it
succeeds, larger gains in gradient's evaluation cost at the price of more costly
function evaluations.

\numsection{Conclusions and Perspectives}\label{conclusions-s}

We have provided an improved provably convergent variant of the trust-region
method using dynamic accuracy and have shown that, when considered in the
context high performance computing and multiprecision arithmetic, this variant
has the potential to bring significant savings in objective function's and
gradient's evaluation cost.

In the deep learning context, computation in reduced accuracy has already
attracted a lot of attention (see \cite{Wangetal18} and references therein,
for instance), but the process, for now, lacks adaptive mechanisms and formal
accuracy guarantees. Our approach can be seen as a step towards providing
them, and the implemenation of our ideas in practical deep learning frameworks
is the object of ongoing research.

Despite the encouraging results reported in this note, the authors are of
course aware that the numerical experiments discussed here are limited in size
and scope and that the suggested conclusions need further assessment. In
particular, it may be of interest to compare inexact trust-region algorithms
with inexact regularization methods \cite{BellGuriMoriToin18}, especially if
not only first-order but also second-order critical points are sought.

\appendix

\appnumsection{Appendix: Complexity Theory for the TR1DA Algorithm}

For the sake of accuracy and completeness, we now provide details of the
first-order worst-case complexity analyis summarized at the end of
Section~\ref{algo-s}. As indicated there, the following development can be
seen as a combination of the arguments proposed by \cite{ConnGoulToin00}
for the convergence theory of trust-region methods with inexact gradients
(pp. 280sq) and dynamic accuracy (pp. 400).

We assume that
\begin{description}
\item[AS.1:] The objective function $f$ is twice continuously differentiable
  in $\Re^n$ and there exist a constant $\kappa_\nabla\geq 0$ such that
  $\|\nabla_x^2f(x)\|\leq \kappa_\nabla$ for all $x \in \Re^n$.
\item[AS.2:] There exists a constant $\kappa_H\geq 0$ such that $\|H_k\| \leq
  \kappa_H$ for all $k\geq 0$.
\item[AS.3] There exists a constant $\kappa_{\rm low}$ such that $f(x)\geq
  \kappa_{\rm low}$ for all $x\in \Re^n$.
\end{description}

\llem{Th8.4.2}{Suppose AS.1 and AS.2 hold. Then, for each $k\geq 0$,
  \beqn{8.4.11}
  |f(x_k+s_k)-m(x_k,s_k)|
  \leq |f_k-f(x_k)| +\kappa_g \|\barg(x_k,\omega_{g,k})\| \Delta_k + \kappa_{H\nabla}\Delta_k^2.
  \eeqn
  for $\kappa_{H\nabla} = 1+\max[\kappa_H, \kappa_\nabla]$.
}

\proof{(See \cite[Theorem 8.4.2]{ConnGoulToin00}.)
  The definition \req{model}, \req{bargk-def}, the mean-value theorem, the Cauchy-Schwarz
  inequality and AS.1 give that, for some $\xi_k$ in the segment $[x_k,x_k+s_k]$,
  \[
  \begin{array}{lcl}
  |f(x_k+s_k)-m(x_k,s_k)|
  & \leq & |f_k-f(x_k)| +|s_k^T(\nabla_x^1f(x_k)-\barg(x_k,\omega_{g,k})| \\*[1.5ex]
  &      & \ms\ms + \half |s_k^T\nabla_x^2f(\xi_k)s_k|
       + \half |s_k^TH_ks_k| \\*[1.5ex]
  & \leq & |f_k-f(x_k)| +\kappa_g \|\barg(x_k,\omega_{g,k})\|\, \|s_k\| + \half(\kappa_H + \kappa_\nabla)\|s_k\|^2
 \end{array}
  \]
  and \req{8.4.11} follows from the the Cauchy-Schwarz inequality and the
  inequality $\|s_k\|\leq \Delta_k$.
} 

\llem{p401}{We have that, for all $k\geq 0$,
  \beqn{ff-bounds}
  \max\left[|f_k-f(x_k)|,|f_k^+ - f(x_k+s_k)|\right] \leq
  \eta_0\left[m(x_k,0)-m(x_k,s_k)\right]
  \eeqn
  and
  \beqn{true-rho} 
  \rho_k \geq \eta_1
  \tim{ implies that }
  \frac{f(x_k)-f(x_k+s_k)}{m(x_k,0)-m(x_k,s_k)} \geq \eta_1 - 2\eta_0 >0.
  \eeqn
}

\proof{(See \cite[p. 401]{ConnGoulToin00}.) The mechanism of the TR1DA
  algorithm ensures that \req{ff-bounds} holds.
  Hence,
  \[
  \frac{\left[f_k-f(x_k)\right]+\left[|f_k^+ -f(x_k+s_k)\right]}{m(x_k,0)-m(x_k,s_k)}
    \leq 2 \eta_0.
  \]
  As a consequence, for iterations where $\rho_k \geq \eta_1$,
  \[
  \rho_k =  \frac{f_k-f_k^+}{m(x_k,0)-m(x_k,s_k)}
  = \frac{f(x_k)-f(x_k+s_k)}{m(x_k,0)-m(x_k,s_k)}
  + \frac{\left[f_k-f(x_k)\right]+\left[|f_k^+-f(x_k+s_k)\right]}
         {m(x_k,0)-m(x_k,s_k)}
  \]
  and \req{true-rho} must hold.
} 

\noindent
This result implies, in particular, that the sequence $\{f(x_k)\}_{k\geq 0}$ is
non-increasing, and the TR1DA algorithm is therefore monotone on the exact
function $f$.

\llem{Th8.4.3}{Suppose AS.1 and AS.2 hold, and that
  $\barg(x_k,\omega_{g,k})\neq 0$. Then
  \beqn{Delta-succ}
  \Delta_k \leq \frac{\|\barg(x_k,\omega_{g,k})\|}{2\kappa_{H\nabla}}
  \Big[ \half(1-\eta_1)-\eta_0-\kappa_g\Big]
  \tim{implies that}
  \Delta_{k+1} \geq \Delta_k.
  \eeqn
}

\proof{(See \cite[Theorem 8.4.3]{ConnGoulToin00}.)
Since \req{alg-consts} implies that $\half(1-\eta_1)-\eta_0-\kappa_g \in (0,1)$
the first part of \req{Delta-succ} then gives that
$\Delta_k < \|\barg(x_k,\omega_{g,k})\|/\kappa_{H\nabla}$. Hence the
inequality $1+\|H_k\|\leq \kappa_{H\nabla}$ and 
\req{cauchy} yield that
\[
m(x_k,0)-m(x_k,s_k)
\geq \half \|\barg(x_k,\omega_{g,k})\|\min\left[\frac{\|\barg(x_k,\omega_{g,k})\|}{\kappa_{H\nabla}},\Delta_k\right]
= \half \|\barg(x_k,\omega_{g,k})\|\,\Delta_k.
\]
As a consequence, we may use \req{rho-def}, the Cauchy-Schwarz
inequality, \req{ff-bounds} (twice), \req{8.4.11}, the inequality
$\kappa_{H\nabla}\geq 1$ and the first part of \req{Delta-succ} to deduce
that, for all $k\geq 0$,
\[
\begin{array}{lcl}
  |\rho_k-1|
  & = & \bigfrac{|f_k^+- m(x_k,s_k)|}{m(x_k,0)-m(x_k,s_k)}\\*[1.5ex]
  & \leq & \bigfrac{|f_k^+-f(x_k+s_k)|+|f(x_k+s_k)-m(x_k,s_k)|}
                   {m(x_k,0)-m(x_k,s_k)}\\*[1.5ex]
  & \leq & 2 \eta_0
     + \bigfrac{\kappa_g\|\barg(x_k,\omega_{g,k})\|\Delta_k+\kappa_{H\nabla}\Delta_k^2}
     {\half \|\barg(x_k,\omega_{g,k})\|\,\Delta_k}\\*[1.5ex]
  & \leq & 2\eta_0 +2\kappa_g+
     2\kappa_{H\nabla}\bigfrac{\Delta_k}{\|\barg(x_k,\omega_{g,k})\|}\\*[1.5ex]
  & \leq & 1-\eta_2.
\end{array}
\]
Thus $\rho_k\geq \eta_2$ and \req{rad-update} ensures the second part of
\req{Delta-succ}.
} 

\llem{Th8.4.4}{Suppose that AS.1 and AS.2 hold.  Then, before termination,
  \beqn{Delta-lower}
  \Delta_k \geq \min\left[\Delta_0,\theta \epsilon\right]
  \tim{ where }
  0< \theta \eqdef
  \frac{\gamma_1\big[\half(1-\eta_1)-\eta_0-\kappa_g\big]}{\kappa_{H\nabla}(1+\kappa_g)}
  <\frac{1}{\kappa_{H\nabla}(1+\kappa_g)}.
  \eeqn
}

\proof{(See \cite[Theorem 8.4.4]{ConnGoulToin00}.)
  Before termination, we must have that
  \beqn{no-term}
  \|\barg(x_k,\omega_{g,k})\|\geq \frac{\epsilon}{1+\kappa_g}.
  \eeqn
  Suppose that iteration $k$ is the first iteration such that
  \beqn{first-D}
  \Delta_{k+1} \leq \frac{\gamma_1\epsilon}{\kappa_{H\nabla}(1+\kappa_g)}
   \big[\half(1-\eta_1)-\eta_0-\kappa_g\big].
  \eeqn
  Then the update \req{rad-update} implies that
  \[
  \Delta_k
  \leq \frac{\epsilon}{\kappa_{H\nabla}(1+\kappa_g)}
       \big[\half(1-\eta_1)-\eta_0-\kappa_g\big]
  \leq \frac{\|\barg(x_k,\omega_{g,k})\|}{\kappa_{H\nabla}}
       \big[\half(1-\eta_1)-\eta_0-\kappa_g\big]
  \]
  where we have used \req{no-term} to deduce the last inequality.
  But this bound and \req{Delta-succ} imply that $\Delta_{k+1}\geq \Delta_k$,
  which is impossible since $\Delta_k$ is reduced at iteration $k$.  Hence no
  $k$ exists such that \req{first-D} holds and the desired conclusion follows.
}

\llem{SvsU-l}{For each $k \geq 0$, define 
  \beqn{SkUk-def}
  \calS_k \eqdef \{ j \in \iibe{0}{k} \mid \rho_j \geq \eta_1 \}
  \tim{ and }
  \calU_k \eqdef \iibe{0}{k} \setminus \calS_k.
  \eeqn
  the index sets of ``successful'' and ``unsuccessful'' iterations,
  respectively. Then
  \beqn{Uk-small}
   k \leq |\calS_k|
   \left(1-\frac{\log\gamma_3}{\log\gamma_2}\right)
   +\frac{1}{|\log\gamma_2|}\log\left(\frac{\Delta_0}{\theta\epsilon}\right).
  \eeqn
}

\proof{Observe that \req{rad-update} implies that
  \[
  \Delta_{j+1} \leq \gamma_3 \Delta_j
  \tim{ for } j \in \calS_k
  \]
  and that
  \[
  \Delta_{j+1}\leq \gamma_2\Delta_j
  \tim{ for } j \in \calU_k.
  \]
  Combining these two inequalities, we obtain from \req{Delta-lower} that
    \[
  \min\big[\Delta_0,\theta\epsilon\Big]\leq \Delta_k \leq
  \Delta_0\gamma_3^{|S_k|}\gamma_2^{|U_k|}
  \]
  Dividing by $\Delta_0$, taking logarithms and recalling that $\gamma_2\in
  (0,1)$, we get
  \[
  |\calU_k| \leq -|\calS_k| \frac{\log\gamma_3}{\log\gamma_2} -
  \frac{1}{\log\gamma_2}\log\left(\frac{\Delta_0}{\theta \epsilon}\right).
  \]
  Hence \req{Uk-small} follows since $k = |\calS_k|+|\calU_k|$.
} 

\lthm{comp-th}{Suppose that AS.1--AS.3 hold. Suppose also that
  $\Delta_0\geq \theta\epsilon$, where $\theta$ is defined in
  \req{Delta-lower}. Then the TR1DA algorithm produces an iterate $x_k$ such
  that $\|\nabla_x^1f(x_k)\| \leq \epsilon$ in at most
  \beqn{succ-its}
  \tau_\calS
  \eqdef \frac{2(f(x_0)-\kappa_{\rm low})(1+\kappa_g)}{(\eta_1-2\eta_0)\theta}
  \cdot \frac{1}{\epsilon^2}
  \eeqn
  successful iterations, at most
  \beqn{total-its}
  \tau_{\rm tot} \eqdef \tau_S \left(1-\frac{\log\gamma_3}{\log\gamma_2}\right)
  +\frac{1}{|\log\gamma_2|}\log\left(\frac{\Delta_0}{\theta\epsilon}\right)
  \eeqn
  iterations in total, at most $\tau_{\rm tot}$ (approximate) evaluations of
  the gradient satisfying \req{bargk-def}, and at most $2\tau_{\rm tot}$
  (approximate) evaluations of the objective function satisfying \req{barf-def}.
}

\proof{
  As in the previous proof, \req{no-term} must hold before termination.
  Using AS.3, \req{SkUk-def}, \req{true-rho}, \req{cauchy}, \req{no-term},
  the assumption that $\Delta_0\geq \theta\epsilon$ and \req{Delta-lower}, we
  obtain that, for an arbitrary $k\geq 0$ before termination,
  \[
  \begin{array}{lcl}
    f(x_0)-\kappa_{\rm low}
    & \geq & \bigsum_{j\in \calS_k}\left[f(x_j)-f(x_{j+1})\right] \\*[1.5ex]
    & \geq & (\eta_1-2\eta_0)\bigsum_{j\in\calS_k}\left[m(x_j,0)-m(x_j,s_j)\right] \\*[1.5ex]
    & \geq & \half(\eta_1-2\eta_0)\bigsum_{j\in\calS_k}\|\barg(x_j,\omega_{g,j})\|
             \min\left[\bigfrac{\|\barg(x_j,\omega_{g,j})\|}{1+\|H_j\|},\Delta_j\right]\\*[1.5ex]
    & \geq & \half|\calS_k|(\eta_1-2\eta_0)\bigfrac{\epsilon}{1+\kappa_g}
             \min\left[\bigfrac{\epsilon}{\kappa_{H\nabla}(1+\kappa_g)},\min\Big[\Delta_0,\theta\epsilon\Big]\right]\\*[1.5ex]
    &   =  & |\calS_k| \,\bigfrac{(\eta_1-2\eta_0)}{2(1+\kappa_g)}
             \min\left[\bigfrac{1}{\kappa_{H\nabla}(1+\kappa_g)},\theta\right]\,\epsilon^2 \\*[1.5ex]
    &   =  & |\calS_k|\,\bigfrac{(\eta_1-2\eta_0)\theta}{2(1+\kappa_g)} \, \epsilon^2
  \end{array}
  \]
  and therefore
  \[
  |\calS_k|
  \leq  \frac{2(f(x_0)-\kappa_{\rm low})(1+\kappa_g)}{(\eta_1-2\eta_0)\theta}
  \cdot \frac{1}{\epsilon^2}
  \eqdef \frac{\tau_\calS}{\epsilon^2}.
  \]
  As a consequence $\|\barg(x_k,\omega_{g,k})\| < \epsilon/(1+\kappa_g)$ after
  at most $\tau_S\epsilon^{-2}$ successful iterations and the algorithm
  terminates. The relation \req{term-grad} then ensures that
  $\|\nabla_x^1f(x_k)\| < \epsilon$, yielding \req{succ-its}. We may now use
  \req{Uk-small} and the mechanism of the algorithm to complete the proof.
}  

\noindent
Given that $\epsilon \in (0,1]$, we immediately note that
\[
\tau_\calS  = \calO(\epsilon^{-2})
\tim{and}
\tau_{\rm tot}  = \calO(\epsilon^{-2}).
\]
Moreover, the proof of Theorem~\ref{comp-th} implies that
these complexity bounds improve from $\calO(\epsilon^{-2})$ to $\calO(\epsilon^{-1})$
if $\epsilon$ is so large or $\Delta_0$ so small to yield $\Delta_0 < \theta\epsilon$.

We conclude this brief complexity theory by noting that the domain in which AS.1 is
assumed to hold can be restricted to the ``tree of iterates''
$\cup_{k\geq0}[x_k,x_k+s_k]$ without altering our results. This can be useful
if an upper bound $\bar{\Delta}$ is imposed on the step's length, in which case the
monotonicty of the algorithm ensures that the tree of iterates remains in the set
\[
\{y \in \Re^n \,|\, y = x+s \tim{with} f(x)\leq f(x_0)
                \tim{and} \|s\|\leq\bar{\Delta} \}.
\]
While it can be difficult to verify AS.1 on the (\emph{a priori} unpredictable) tree
of iterates, verifying it on the above set is much easier.
\end{document}